\documentclass[twocolumn]{autart}
\usepackage{color}
\usepackage{indentfirst,comment} 
\usepackage{graphicx}
\usepackage{amsmath}
\usepackage{amsfonts}
\usepackage{mathtools}
\usepackage{algorithm}
\usepackage{algpseudocode}
\usepackage{tikz}
\usepackage{adjustbox}
\usepackage{url}
\usepackage{enumitem}
\usepackage[normalem]{ulem}
\usepackage{cancel}

\newcommand{\rev}[1]{\textcolor{black}{#1}}
\newcommand{\revv}[1]{\textcolor{black}{#1}}
\newtheorem{defi}{Definition}
\newtheorem{assmp}{Assumption}
\newtheorem{rmk}{Remark}

\newcommand{\bignorm}[1]{\left\lVert#1\right\rVert}
\newcommand{\software}[1]{{\tt#1}}
\DeclarePairedDelimiter{\ceil}{\lceil}{\rceil}

\algnewcommand{\Initialize}[1]{
  \State \textbf{Initialize:}
  \Statex \hspace*{\algorithmicindent}\parbox[t]{.95\linewidth}{\raggedright #1}
}

\algnewcommand{\Input}[1]{
  \State \textbf{Input:}
  \Statex \hspace*{\algorithmicindent}\parbox[t]{.95\linewidth}{\raggedright #1}
}

\begin{document}
\begin{frontmatter}
\title{An Efficient MPC Algorithm For Switched Systems with Minimum Dwell Time Constraints}

\author[FZU]{Yutao Chen}\ead{yutao.chen@fzu.edu.cn} 
\author[TUE]{Mircea Lazar}\ead{m.lazar@tue.nl}

\address[FZU]{College of Electrical Engineering and Automation, Fuzhou University, China.}
\address[TUE]{Department of Electrical Engineering, Eindhoven University of Technology, the Netherlands.}



\begin{abstract}
This paper presents an efficient suboptimal model predictive control (MPC) algorithm for nonlinear switched systems subject to minimum dwell time constraints (MTC). While MTC are required for most physical systems due to stability, power and mechanical restrictions, MPC optimization problems with MTC are challenging to solve. To efficiently solve such problems, the on-line MPC optimization problem is decomposed into a sequence of simpler problems, which include two nonlinear programs (NLP) and a rounding step, as typically done in mixed-integer optimal control (MIOC). Unlike the classical approach that embeds MTC in a mixed-integer linear program (MILP) with combinatorial constraints in the rounding step, our proposal is to embed the MTC in one of the NLPs using move blocking. Such a formulation can speedup on-line computations by employing recent move blocking algorithms for NLP problems and by using a simple sum-up-rounding (SUR) method for the rounding step. An explicit upper bound of the integer approximation error for the rounding step is given. In addition, a combined shrinking and receding horizon strategy is developed to satisfy closed-loop MTC. Recursive feasibility is proven using a $l$-step control invariant ($l$-CI) set, where $l$ is the minimum dwell time step length. An algorithm to compute $l$-CI sets for switched linear systems off-line is also presented. Numerical studies \revv{show significant speed-up and comparable control performance of the proposed MPC algorithm against the classical approach, though at the cost of sub-optimal solutions.}
\end{abstract}

\end{frontmatter}

\section{Introduction}
Switched systems are a special class of hybrid systems that consist of a number of modes out of which only one is active at a given time. The switch from one mode to another can be triggered by an external control input or by certain internal conditions \cite{zhu2015optimal}. Optimal control of switched systems formulates an optimal control problem (OCP) where a cost function is minimized to find the optimal switching strategy as well as the state and continuous input trajectory. The switching strategy includes the sequence of switching modes and the sequence of time instances at which switching occurs. Applications of optimal control of switched systems can be found in mode scheduling for automobiles, valve control for chemical processes, pesticide scheduling in agriculture, and many other applications, see, for example, \cite{zhu2015optimal}, \cite{alamir2004solving}, \cite{ali2018hybrid} and references therein.

When the OCP is solved on-line repeatedly using a finite prediction horizon and the latest measured or estimated data, it gives rise to model predictive control (MPC), which is a widely used advanced control technique. The key factors for the success of MPC are  efficient algorithms for solving the OCP in real-time and conditions for guaranteeing recursive feasibility and closed-loop stability. However, the switching dynamics makes the on-line optimization challenging, especially when minimum dwell time constraints (MTC) exist. MTC limit the minimum time of any active mode before switching to another. MTC are required in many real-world applications, e.g. shifting the gear of a vehicle needs a noteworthy amount of time. MTC are also important when considering that the number of switches cannot be infinite and the switching cost cannot be ignored.

\subsection{Relevant Work}
As far as the authors know, using MPC for switched systems to recursively determine the switching sequence subject to MTC has only been recently studied in \cite{burger2019design}. The main idea proposed therein is as follows. First, the on-line optimization problem has been transcribed into a mixed-integer nonlinear program (MINLP) using direct methods. Then, a decomposition method developed in the mixed-integer optimal control (MIOC) field has been adopted, where the MINLP is decomposed into several simpler problems including two nonlinear programs (NLP) and a mixed-integer linear program (MILP). The MTC are removed from the NLPs and are formulated as combinatorial constraints in the corresponding MILP. To satisfy the closed-loop MTC, at each sampling instant, the first mode of the open-loop mode sequence is fixed to the previous active mode if MTC are active while the terminal mode of the mode sequence is not constrained. 

It is worth to mention that a real-time MPC algorithm for switched nonlinear systems has been proposed in \cite{katayama2020moving}, but without considering MTC or other types of constraints. Also, \rev{many efforts have} been dedicated to develop efficient numerical solvers for MPC of switched systems without MTC, starting with the seminal work \cite{bemporad1999control} for mixed-logical dynamical systems and \rev{followed by} the recent advances in branch and bound methods for solving MINLP, see \cite{bemporad2018numerically}, \cite{hespanhol2019structure} and the references therein.  
Alternative numerical methods that could be used to solve on-line optimization problems within MPC for switched systems have been developed for decades in the field of optimal control of switched systems. These include studies assuming a fixed switching sequence \cite{xu2004optimal}, \cite{egerstedt2006transition},  hybrid minimum principle based methods \cite{alamir2004solving}, \cite{shaikh2007hybrid}, \cite{gonzalez2010descent} and  mode insertion algorithms \cite{wardi2015switched}. Efficient optimal control algorithms using direct methods have been developed in \cite{bengea2005optimal}, \cite{sager2009reformulations}. Both these methods have adopted a similar idea of first relaxing \rev{the integrality constraints} and then rounding. To deal with MTC, the authors of \cite{ali2014optimal} have adopted mode insertion techniques. This approach has been further extended to consider different MTC for different modes in \cite{ali2018hybrid}, given a fixed mode sequence. A method based on approximate dynamic programming (ADP) has been proposed in \cite{heydari2017optimal} for unconstrained systems. General constrained systems with MTC have been tackled in \cite{burger2015dynamic} using dynamic programming (DP) which suffers from the curse of dimensions. In MIOC, a modified rounding strategy has been proposed to reduce the complexity of the rounding step at the cost of loss of optimality \cite{zeile2020mixed}.

It should be noted that another category of studies on MPC for switched systems have focused on developing stabilizing MPC controllers. Interested readers can refer to \cite{mhaskar2005predictive}, \cite{lazar2006stabilizing}, \cite{muller2012model}, \cite{zhang2016switched} and references therein. The main difference from this category to the optimal control of switched systems is that, the switching sequence in the stabilizing MPC controllers, either known or unknown a priori, is not a decision variable. Another difference is that these papers focus on exploiting knowledge of the dwell time, e.g. the average dwell time specified in \cite{muller2012model}, for ensuring stability, without addressing the question of how to generate feasible switching sequences that comply with MTC. \rev{Such a dwell time strategy has also been adopted in \cite{franze2017command} using the command governor structure, where a prescribed dwell time is computed offline to ensure stability and constraint fulfillment during the switching transient phase. In contrast, this paper considers MTC as a hard constraint when deciding the switching between modes on-line, which is assumed instantaneous, i.e., without transient dynamics.}

\subsection{Contribution}
In this paper, we develop an efficient suboptimal MPC algorithm for computing the switching sequence, the state and control trajectory for constrained nonlinear switched systems subject to MTC. To solve the OCP on-line, we propose a variant of the MIOC decomposition method \cite{burger2019design}, \cite{kirches2020approximation}, which embeds the MTC in the NLP instead of the corresponding MILP. In particular, the integer control variables which represent the switching sequence are relaxed and ``blocked'' (i.e., set constant) over multiple discretization intervals using move blocking techniques. The relaxed solutions of the NLP are then rounded using a simple sum-up-rounding (SUR) method with bounded integer approximation error. 
The computational effort to solve the NLP can be reduced by employing efficient move blocking numerical algorithms \cite{chen2019efficient}. In addition, the SUR step has negligible computational burden when compared to solving an MILP with combinatorial constraints. 

A second main contribution of this paper is a combined, shrinking horizon and receding horizon terminal set strategy for ensuring recursive feasibility of the developed MPC algorithm. 
We show that this combination is recursively feasible when a $l$-step control invariant ($l$-CI) \cite{lazar2013stability} terminal set is employed. For switched linear systems, an algorithm is developed to explicitly compute a $l-$step switch-robust CI ($l$-SRCI) set, which is a practical class of $l$-CI sets.
Finally, we implement the proposed method in \software{MATMPC}, an open-source nonlinear MPC tool based on \software{MATLAB} that supports tailored move blocking algorithms \cite{chen2019matmpc}. 
\revv{We show that the developed MPC algorithm achieves significant speed-up over the existing MILP based MPC algorithms for two MIOC benchmark examples, at the expense of sub-optimal solutions (i.e., due to the reduced degrees of freedom regarding switching). However, closed-loop simulations show comparable state trajectories and objective function values despite sub-optimal solutions.}


\section{Problem Description and Preliminaries}\label{sec2}
In this section we first present definitions, assumptions and the considered OCP formulation. Then we introduce recent results from the MIOC field on how to solve the OCP.

\subsection{Definitions, assumptions and problem formulation}\label{sec2-subsec1}
Consider switching mode continuous time nonlinear dynamical systems 
\begin{equation}
    \dot{x}(t) = f_{q(t)}(x(t),u(t)),\quad x(t_0)=x_0,\quad \revv{t\in\mathbb{R}_+},
\end{equation}
where $x(t)\in\mathbb{R}^{n_x}, u(t)\in\mathbb{R}^{n_u}$ are 
continuous-time state and input trajectories and $f_{q(t)}:\mathbb{R}^{n_x}\times \mathbb{R}^{n_u}\rightarrow \mathbb{R}^{n_x}$ is $\mathcal{C}^2$. Above $\mathbb{R}_+$ denotes the set of non-negative real numbers. We adopt the following definitions.
\begin{itemize}
    \item The active mode function $q(t):[t_0,t_f]\rightarrow \mathbb{Q}$ is piece-wise constant and left-continuous in a given and fixed time period $[t_0,t_f]\subset\mathbb{R}_+$, where $\mathbb{Q}\coloneqq\{1,2,\ldots,Q\}, \,Q\in\mathbb{N}$ is a finite set of available modes.
    \item Let $K\in\mathbb{N}$ denote the number of mode switches during the time interval $[t_0,t_f]$.
    \item Define the set of switching sequences by $\mathcal{Q}=\{(q_1,q_2,\ldots,q_{K+1})^\top, q_i\in\mathbb{Q}, \,\forall i=1,\ldots,K+1\}$.
    \item Define the set of switching time instants by $\mathcal{T}\coloneqq\{\tau\in \rev{[t_0,t_f]^{K+2}}:\tau_i\leq \tau_{i+1},\, \rev{i=0},\ldots,K\}$, where \rev{$\tau_0=t_0$}, $\tau_{K+1}=t_f$ and $\tau\coloneqq(\tau_0,\tau_1,\tau_2,\ldots,\tau_{K+1})^\top$. 
    \item Define the set of switching laws by
    \begin{equation}
        \Sigma = \{(\rev{q(\cdot)},\tau): \rev{q(\cdot)}\in\mathcal{Q}, \tau\in \mathcal{T}\}.
    \end{equation}
    Let $\rev{\sigma(\cdot)}=(\rev{q(\cdot)},\tau)\in\Sigma$ denote a switching law.
\end{itemize}
We make the following assumptions.
\begin{assmp}
There is no state jump at switching time instances, i.e., $x_{q_{i+1}}(\tau_{i+1})=\lim_{t\rightarrow \tau_{i+1}}x_{q_{i}}(t)$, for all $\rev{i=0},\ldots,K-1.$
\end{assmp}

\begin{assmp}\label{assmp2}
There is no autonomous switching, i.e. the mode switching is controlled by $q(t)$, which is a signal to be determined.
\end{assmp}


\revv{Let $\tilde{x}_0$ be given and consider the following OCP defined over a time interval $[t_0,\, t_f]\subset\mathbb{R}_+$:}
    \begin{subequations}\label{OCP}
		\begin{align}
			\min_{x(t),u(t),\sigma(t)} \quad & \psi(x(t_f)) +\sum_{i=0}^K \int_{\tau_i}^{\tau_{i+1}}L_{q(t)}(x(t),u(t))\text{d}t\\
			s.t. \quad & \dot{x}(t) = f_{q(t)}(x(t),u(t)),\,x(t_0)=\tilde{x}_0,\\
			& x(t)\in\mathcal{X},\,\rev{\forall t\in[t_0,t_f],}\\
			& u(t)\in\mathcal{U},\,\rev{\forall t\in[t_0,t_f],}\\
			& \sigma(t)\in\Sigma,\,\rev{\forall t\in[t_0,t_f],}\\
			& g(x(t),u(t))\leq 0,\,\rev{\forall t\in[t_0,t_f],}\label{path con}\\
			& x(t_f)\in\mathcal{X}_f,\\
			& 0<\underline{\Delta\tau}\leq \tau_{i+1}-\tau_i\leq t_f-t_0, \nonumber\\
			&\quad\forall i=0,\ldots,K,\label{mini dwell}
		\end{align}
	\end{subequations}
	where $\mathcal{X}$ is a closed set and $\mathcal{U}$ is compact set. The running cost $L_{q(t)}:\mathbb{R}^{n_x}\times \mathbb{R}^{n_u}\rightarrow \mathbb{R}$ and the terminal cost $\psi:\mathbb{R}^{n_x}\rightarrow \mathbb{R}$ are assumed to be $\mathcal{C}^2$ functions.  The initial condition $\tilde{x}_0$ is the measured or estimated state. The function $g: \mathbb{R}^{n_x}\times \mathbb{R}^{n_u}\rightarrow \mathbb{R}^c$ is assumed to be $\mathcal{C}^2$. The set $\mathcal{X}_f\subseteq \mathcal{X}$ defines the terminal constraint. MTC are defined by \eqref{mini dwell} for a given $\underline{\Delta\tau}$. Note that $g$ can also depend on the mode $q$, which results in a \revv{mathematical program with vanishing constraints (MPVC)} \cite{izmailov2009mathematical}. We do not consider this situation in this paper.

\subsection{Mixed-integer optimal control}\label{sec2-subsec2}
In this paper we focus on efficiently solving problem \eqref{OCP} via finite dimensional approximations. To this end, one of the promising approaches is to approximate \eqref{OCP} by a MIOC problem, in which the \rev{switching sequence is transformed into a set of binary functions using outer convexification} \cite{sager2009reformulations}.  Given binary functions $b^j(t):[t_0,t_f]\rightarrow\{0,1\}, j=1,\ldots,Q$, the MIOC problem is formulated as 
\begin{subequations}\label{MIOCP}
\begin{align}
     \min_{x(t),u(t),b(t)} \quad & \psi(x(t_f)) + \int_{t_0}^{t_f}\sum_{j=1}^Q b^j(t)\rev{L_j}(x(t),u(t))\text{d}t, \label{obj}\\
    s.t. \quad & \dot{x}(t) = \sum_{j=1}^Q b^j(t)f_{j}(x(t),u(t)),\nonumber\\ &x(t_0)=\tilde{x}_0, \label{nonlinear dynamics}\\
        & x(t)\in\mathcal{X},\,\forall t\in[t_0,t_f],\label{state cons}\\
        & u(t)\in\mathcal{U},\,\forall t\in[t_0,t_f],\\
        & x(t_f)\in\mathcal{X}_f,\,\forall t\in[t_0,t_f],\\
        & g(x(t),u(t))\leq 0,\,\forall t\in[t_0,t_f],\label{path con miocp}\\
        & \sum_{j=1}^Q b^j(t)=1, \,\forall t\in[t_0,t_f], \label{SOS1}
\end{align}
\end{subequations}
where \revv{the binary functions can be stacked into a vector-valued mapping}
\begin{equation}\label{definition b}
    \mathbf{b}\revv{:[t_0,t_f]\rightarrow \{0,1\}^Q},\quad \mathbf{b}(t):=(b^1(t),\ldots,b^Q(t))^\top.
\end{equation}
The solutions to problem \eqref{MIOCP} and \eqref{OCP} are bijective without constraints \eqref{mini dwell} \cite{sager2009reformulations}. The special ordered set of type~1 (SOS1) \eqref{SOS1} guarantees that only one mode is active at any time \cite{sager2009reformulations}.

\rev{This paper follows the decomposition method proposed in \cite{sager2009reformulations}, \cite{burger2019design} that breaks down \eqref{MIOCP} into a sequence of simpler problems.} The procedure consists of three steps\footnote{In \cite{robuschi2021multiphase}, the three steps are repeated multiple times to improve feasibility of the optimization problem.}: 
\begin{enumerate}[label=(\roman*)]
    \item NLP \rev{\#1}: discretize \eqref{MIOCP} and replace the binary mapping \revv{$\mathbf{b}(t)$} with \revv{the real-valued mapping $\hat{\mathbf{b}}:[t_0,t_f]\rightarrow[0,1]^Q$};
    \item Combinatorial integral approximation (CIA): \revv{given $\hat{\mathbf{b}}(t)$}, solve an MILP to satisfy MTC constraints \eqref{mini dwell}, obtaining $\mathbf{b}(t)$;
    \item NLP \rev{\#2}: solve NLP \#1 with fixed binary variables $\mathbf{b}(t)$ obtained from the CIA step (ii), \revv{obtaining $x(t),u(t)$}.
\end{enumerate}

\section{MPC algorithm for switched systems with MTC using move blocking}\label{sec3}
In this section, we develop an MPC algorithm by embedding MTC into the NLP \#1 \rev{hence avoiding formulating and solving an MILP in the CIA step}. This is achieved by imposing move blocking on the \revv{real-valued mapping $\hat{\mathbf{b}}(t)$}. As a consequence, the binary mapping \revv{$\mathbf{b}(t)$} automatically satisfies the MTC.

\subsection{The NLP \#1 with move blocking}
To numerically solve \eqref{MIOCP}, we adopt the direct multiple shooting method \cite{bock1984multiple}. The time domain $[t_0,t_f]$ is discretized into $N$ intervals, characterized by equidistant grid points
\begin{equation}\label{grid}
    t_0 < t_1 < \ldots < t_{N-1} < t_N = t_f,
\end{equation}
with $\Delta t=\frac{t_f-t_0}{N}$. In each interval, the control and binary functions $(u(t),\mathbf{b}(t))$ are assumed to be constant, hence they can only change values at grid points. Define the minimum dwell time interval length as
\begin{equation}\label{compute l}
    l := \ceil{\frac{\underline{\Delta \tau}}{\Delta t}}
\end{equation}
where \rev{$\ceil{}$ is the ceil function}. To meet the MTC, \revv{we introduce move blocking for $l$ consecutive real-valued vectors $\hat{\mathbf{b}}_k$ for $k=0,l,\ldots,(M-1)l$, as follows:}
\revv{
\begin{equation}\label{input block}
    \hat{\mathbf{b}}_{k+i} = \hat{\mathbf{b}}_{k},\quad \forall i=1,\ldots, l-1,
\end{equation}}
\revv{where $\hat{\mathbf{b}}_k=\hat{\mathbf{b}}(t_k)\in[0,1]^Q$ is the discretized real-valued mapping at grid point $t_k$.}
Note that $M<K\in\mathbb{N}$ is the number of ``blocks'' (with each block consisting of $l$ equal elements, i.e., $\{\hat{\mathbf{b}}_{k+i}\}_{i=0,\ldots,l-1}$) within the prediction horizon, which implicitly specifies the allowed number of switches. Note that we can always find an appropriate tuning parameter $N$ and an approximate upper bound on $\underline{\Delta \tau}$ such that $Ml=N$. The constraints \eqref{input block} enforce $\hat{\mathbf{b}}_k$ to be constant over $l$ consecutive intervals. 
\revv{To reduce the number of decision variables, introduce a real-valued mapping $\hat{\mathbf{p}}:[t_0,t_f]\rightarrow [0,1]^Q$, such that
\begin{align}\label{definition p}
    \hat{\mathbf{p}}_m\coloneqq(\hat{p}_m^1,\ldots,\hat{p}_m^Q)^\top\in[0,1]^Q, \,m=0,1,\ldots,M-1, 
\end{align}
represents the $M$ decision variables $\{\hat{\mathbf{b}}_{k}\}_{i=0,l,\ldots,(M-1)l}$ in \eqref{input block}.} As a result, we obtain the NLP as:
\begin{subequations}\label{NLP-blocked}
\begin{align}
     \min_{\mathbf{x},\mathbf{u},\revv{\hat{\mathbf{p}}}} \quad & \psi(x_N) + \sum_{k=0}^{N-1} \sum_{j=1}^Q \revv{\hat{p}_{m}^{j}} L_j(x_k,u_k) \\
    s.t. \quad & x_{k+1} = \phi(x_k,u_k,\hat{\mathbf{p}}_m), \\
               & \quad \quad k=0,1,\ldots,N-1,\nonumber\\
               & \quad \quad m=0,1,\ldots,M-1,\nonumber\\
        & x_0 = \tilde{x}_0,\\
        & x_k\in\mathcal{X},\,\rev{k=0,1,\ldots,N-1,}\\
        & u_k\in\mathcal{U},\,\rev{k=0,1,\ldots,N-1,}\\
        & x_{N}\in\mathcal{X}_f,\\
        & g(x_k,u_k)\leq 0,\,\rev{k=0,1,\ldots,N-1,}\\
        & \revv{\sum_{j=1}^Q \hat{{p}}_{m}^{j}=1,\,m=0,1,\ldots,M-1,}
\end{align}
\end{subequations}
where 
\begin{align*}
    &\mathbf{x} = (x_0^\top,x_1^\top,\ldots,x_N^\top)^\top\revv{\in\mathbb{R}^{n_x(N+1)}},\\
    &\mathbf{u} = (u_0^\top,u_1^\top,\ldots,u_{N-1}^\top)^\top\revv{\in\mathbb{R}^{n_uN}},\\
    &\revv{\hat{\mathbf{p}} = (\hat{\mathbf{p}}_0^\top,\hat{\mathbf{p}}_1^\top,\ldots,\hat{\mathbf{p}}_{M-1}^\top)^\top\in\mathbb{R}^{MQ},}
\end{align*}
and \revv{$\hat{p}_{m}^{j}$ is $j-$th element of $\hat{\mathbf{p}}_m$}. \revv{The dependency of $\mathbf{x}, \mathbf{u}$ and $\mathbf{p}$ on the current discrete-time instant was omitted above for simplifying the notation.} The function $\phi: \mathbb{R}^{n_x}\times \mathbb{R}^{n_u} \times [0,1]^{Q}\rev{\rightarrow \mathbb{R}^{n_x}}$ is \rev{a numerical integration operator on the nonlinear dynamics \eqref{nonlinear dynamics} using methods such as Euler's and Runge-Kutta}.

\subsection{Integer approximation error}
After solving \eqref{NLP-blocked}, a simple SUR step can be employed to obtain the corresponding binary mapping \revv{$\mathbf{p}:[t_0,t_f]\rightarrow \{0,1\}^Q$ from $\hat{\mathbf{p}}$} \cite{sager2009reformulations}. However, the block constraints \eqref{input block} essentially change the discretization interval length from \revv{$\Delta t$ for $\hat{\mathbf{b}}$, to $l\Delta t$ for $\hat{\mathbf{p}}$}. This leads to minor changes to the SUR scheme as well as the upper bound of the integer approximation error. For $m=0,1,\ldots,M$, the SUR scheme computes

\rev{
\begin{subequations}\label{SUR}
\begin{align}
    &s_{m}^{j} = \left(\sum_{r=0}^m \hat{p}_{r}^{j} - \sum_{r=0}^{m-1}p_{r}^{j}\right)l\Delta t, \,j=1,\ldots,Q,\\
    &p_{m}^{j} = \begin{cases}
                1 \quad \text{if } s_{m}^{j}\geq s_{m}^{d}\,\left(\forall d\neq j\right) \& \left(j<d,\,\forall d:s_{m}^{j}=s_{m}^{d}\right)\\
                0 \quad \text{otherwise},
              \end{cases}
\end{align}
\end{subequations}
where $m,r$ are the indices for blocks and $j,d$ are the indices for vector elements, \revv{i.e. $p_m^j$ is the $j$-th element of $\mathbf{p}_m$.} The auxiliary vector $s$ stores the sum-up difference between the binary and the relaxed variable.} We have the following result.

\begin{prop}\label{thm1}
If $\hat{\mathbf{p}}(t):[t_0,t_f]\rightarrow [0,1]^Q$ is measurable \revv{and essentially bounded}, and $\sum_{j=1}^Q\hat{p}^j(t)=1$ holds, then the function $\mathbf{p}(t):[t_0,t_f]\rightarrow \{0,1\}^Q$ converted from \eqref{SUR} using zero-order hold \revv{on the grid \eqref{grid}, using the block representation \eqref{input block} and \eqref{definition p}, with $Ml=N$}, for $Q\geq 2$, satisfies
\revv{
\begin{equation}\label{CIA error}
\begin{aligned}
    \bignorm{\int_{t_0}^{t} \hat{\mathbf{p}}(\tau)-\mathbf{p}(\tau)\,\text{d}\tau}_{\infty}&\leq l\Delta t\sum_{c=2}^C\frac{1}{c}, \\
\end{aligned}
\end{equation}
}
\revv{where $C=\min\{Q,N+1\}$,} and $\mathbf{p}(t)$ satisfies $\sum_{j=1}^Q p^j(t)=1$.
\end{prop}

\revv{
\begin{pf}
The upper bound for the non-blocked integral approximation comes from Theorem 6.1 of \cite{kirches2020approximation}, where a tightest bound possible for the SUR rounding scheme \eqref{SUR} is proven, written as
\begin{equation*}
\begin{aligned}
    \bignorm{\int_{t_0}^{t} \hat{\mathbf{b}}(\tau)-\mathbf{b}(\tau)\,\text{d}\tau}_{\infty}&\leq \Delta t\sum_{c=2}^C\frac{1}{c}, \\
\end{aligned}
\end{equation*}
where $C=\min\{Q,N+1\}$, and $\mathbf{b}(t), \mathbf{\hat{b}}(t)$ are \revv{the non-blocked corresponding binary and real-valued mappings}. Imposing the blocked decision variables representation \eqref{definition p} with $Ml=N$ extends the grid interval length from $\Delta t$ to $l\Delta t$ for the blocked mappings $\mathbf{p}(t), \mathbf{\hat{p}}(t)$, which leads to the inequality \eqref{CIA error}. \hfill{$\qed$}
\end{pf}
}

\revv{Note that, according to \eqref{compute l}, we have $l\Delta t\leq \underline{\Delta \tau}+\Delta t$.} Therefore, Proposition \ref{thm1} shows that the minimum dwell time $\underline{\Delta \tau}$ enters linearly into the upper bound of the integer approximation error \eqref{CIA error}. Such an upper bound is fixed once the MTC are specified, and it is not relevant to the discretization interval length. Hence, the upper bound of the integer approximation error cannot be made arbitrarily small for the proposed algorithm. Similarly, error upper bounds can also be obtained for the state trajectory $x(t)$, the objective and the path constraint (see Corollary 6 and 8 in \cite{sager2012integer}). It should be noted that another upper bound can be obtained by solving the NLP without move blocking and then applying a modified SUR strategy as in \cite{zeile2020mixed}. Nevertheless, the developed upper bound therein still contains $\underline{\Delta \tau}$ linearly and hence, it cannot be made arbitrarily small either.

A direct consequence of Proposition \ref{thm1} is that the NLP~\#2 may be infeasible. Given that the original \rev{discretized} MIOC problem is feasible, feasibility of the NLP \#1 \eqref{NLP-blocked} is guaranteed due to its larger feasible set. If MTC are not present, the NLP \#2 can be rendered feasible by choosing a sufficiently fine discretization grid with a sufficiently small integer approximation error \cite{sager2012integer}. However, this is no longer the case if MTC are present. In this paper, since the MTC is considered a hard constraint that can be imposed by users, \revv{the minimum dwell time $\underline{\Delta \tau}$ must be small enough to ensure feasibility. The notion of $\delta$-feasibility is defined next.}
\revv{
\begin{defi}
Let $(\hat{x}(t),u^*(t))$ be a feasible state and input trajectory from the solution of \eqref{MIOCP} with $\hat{\mathbf{p}}(t)$ the optimal solution of \eqref{NLP-blocked} and $J(\hat{x},u^*)$ the corresponding cost function, and $x(t)$ a state trajectory for the same input $u^*(t)$ with $\mathbf{p}(t)$ computed from \eqref{SUR} and $J(x,u^*)$ the corresponding cost function. The trajectory $(x(t),u^*(t))$ and problem \eqref{NLP-blocked} is said to be $\delta$-feasible if  $\exists\delta_1, \delta_2, \delta_3\geq 0$ such that
\begin{subequations}\label{defi1 equa}
\begin{align}
    &\bignorm{x(t)-\hat{x}(t)}_{\infty}\leq \delta_1,\label{state CIA}\\
    &|J(x,u^*)-J(\hat{x},u^*)|\leq \delta_2,\\
    &\bignorm{g(x,u^*)-g(\hat{x},u^*)}_{\infty}\leq \delta_3.
\end{align}
\end{subequations}
\end{defi}
We make the following assumption.
\begin{assmp}\label{assmp:feasible dwell time}
The minimum dwell time $\underline{\Delta \tau}$ satisfies the upper bounded 
\begin{equation}
    \underline{\Delta \tau} \leq \frac{\epsilon}{\sum_{c=2}^C\frac{1}{c}}-\Delta t,
\end{equation}
for a given $\epsilon>0$ such that \eqref{CIA error} holds, and the trajectory $(x(t),u^*(t))$ is $\delta$-feasible.
\end{assmp}
}

Assumption \ref{assmp:feasible dwell time} implies that the trajectory $(\hat{x}(t),u^*(t))$ from problem \eqref{NLP-blocked} is not only in the interior of the feasible region, but it also lies at a non-zero distance characterized by constants $\delta_1$ and $\delta_3$, from the boundary of the feasible region, such that the trajectory  $(x(t),u^*(t))$ obtained after the SUR step \eqref{SUR} is feasible.

An immediate advantage of using the proposed move blocking strategy for the NLP \#1 is that applying \eqref{SUR} is computationally much cheaper than solving an MILP for the CIA step. As a result, solving the MIOC problem \eqref{MIOCP} requires only solving two NLPs. In addition, efficient move blocking algorithms by exploiting the move blocking structure to accelerate solving \eqref{NLP-blocked} have been reported, e.g. as in \cite{chen2019efficient}, \cite{son2020move}.


\begin{rmk}
	If system dynamics \eqref{state cons}  is autonomous, i.e. the continuous control input $u(t)$ is absent, there is no need to formulate and solve NLP \#2. After the SUR step \eqref{SUR}, the state trajectory can be obtained by simulating \eqref{state cons} using the calculated binary mode variable.
\end{rmk}

\subsection{The MPC algorithm}
In MPC, problem \eqref{NLP-blocked} must be solved repeatedly on-line in a receding horizon fashion. In particular, the solution of \eqref{NLP-blocked} only guarantees open-loop but not closed-loop MTC fulfillment. To solve this problem, in \cite{burger2019design}, the first mode of the open-loop mode sequence, i.e. $\mathbf{\hat{b}}_0$, is fixed to the previous active mode if the closed-loop MTC are still active. However, as the prediction horizon recedes, the terminal mode of the mode sequence, i.e. $\mathbf{\hat{b}}_{N-1}$ is not constrained and hence, it may violate the MTC. In this work, we propose a combination of shrinking horizon and receding horizon strategies to ensure closed-loop MTC fulfillment with recursive feasibility guarantees.


\revv{Consider in a closed-loop control scenario solving \eqref{NLP-blocked} at time $t_i$ for the sampling instant $i$. Given the current active mode $b_{act}$ and its starting time $t_a\leq t_i$, define the active time duration for $b_{act}$ as $t_{act}\coloneqq t_i-t_a\in\mathbb{R}_{+}$.}
\revv{Define 
\begin{equation}\label{compute-h}
    h \coloneqq \begin{cases}
        \ceil{\frac{\underline{\Delta \tau}-t_{act}}{\Delta t}},\,\text{if} \underline{\Delta \tau}\geq t_{act},\\
        0,\, \text{if} \underline{\Delta \tau}< t_{act},
    \end{cases}
\end{equation}
}
as the number of discretization intervals during which the current mode must \revv{remain the same}. Construct the following move blocking structure and block representation:
\begin{subequations}\label{MPC input block}
\begin{align}
    &\hat{\mathbf{p}}_0 \gets \hat{\mathbf{b}}_0=\ldots=\hat{\mathbf{b}}_{h-1},\label{first control fixed}\\
    &\hat{\mathbf{p}}_1 \gets \hat{\mathbf{b}}_h=\ldots=\revv{\hat{\mathbf{b}}_{h+l-1}}, \label{remaining control blocked}\\
    &\vdots\nonumber\\
    &\hat{\mathbf{p}}_{M-1} \gets \revv{\hat{\mathbf{b}}_{h+(M-2)l}=\ldots=\hat{\mathbf{b}}_{h+(M-1)l-1}}.
\end{align}
\end{subequations}
\revv{Above, the first blocked mode $\hat{\mathbf{p}}_{0}$ represents $h$ consecutive real-valued modes and} is fixed by the current active mode $b_{act}$, \revv{while the remaining blocked modes represent $l$ consecutive real-valued modes}. At sampling instant $i$, define an optimization problem $P_i$ for given parameters $(h, b_{act})$ as:
\begin{subequations}\label{NLP-parametric}
\begin{align}
P_i(h, b_{act}):&\nonumber\\
\min_{\mathbf{x}_i,\mathbf{u}_i,\hat{\mathbf{p}}_i} \quad & \psi(x_{N|i}) + \sum_{k=0}^{N-1-h} \sum_{j=1}^Q \hat{p}_{m|i}^j L_j(x_{k|i},u_{k|i}) \\
    s.t. \quad & \hat{\mathbf{p}}_{0|i} = b_{act},\\
        & x_{k+1|i} = \phi(x_{k|i},u_{k|i},\hat{\mathbf{p}}_{m|i}),\\
            & \quad \quad k=0,1,\ldots,N-1,\nonumber\\
            & \quad \quad m=0,1,\ldots,M-1,\nonumber\\
        & x_{0|i} = \tilde{x}_{0|i},\\
        & x_{k|i}\in\mathcal{X},\,k=0,1,\ldots,N-1,\\
        & u_{k|i}\in\mathcal{U},\,k=0,1,\ldots,N-1,\\
        & x_{N|i}\in\mathcal{X}_f,\\
        & g(x_{k|i},u_{k|i})\leq 0,k=0,1,\ldots,N-1,\\
        & \sum_{j=1}^Q \hat{p}_{m|i}^j=1,\,m=0,1,\ldots,M-1,
\end{align}
\end{subequations}
where 
\revv{
\begin{align*}
    &\mathbf{x}_i = (x_{0|i}^\top,x_{1|i}^\top,\ldots,x_{N|i}^\top)^\top\in\mathbb{R}^{n_x(N+1)},\\
    &\mathbf{u}_i = (u_{0|i}^\top,u_{1|i}^\top,\ldots,u_{N-1|i})^\top\in\mathbb{R}^{n_uN},\\
    &\hat{\mathbf{p}}_i = (\hat{\mathbf{p}}_{0|i}^\top,\hat{\mathbf{p}}_{1|i}^\top,\ldots,\hat{\mathbf{p}}_{M-1|i}^\top)^\top\in\mathbb{R}^{MQ}.
\end{align*}
}
The variable $x_{k|i}$ is the same as $x_k$ in \eqref{NLP-blocked} with an explicit mark on the sampling instant $i$, and $\tilde{x}_{0|i}$ denotes the measured or estimated state at the sampling instant $i$. Algorithm~\ref{alg:mpc} summarizes the proposed MPC scheme, where problem \eqref{NLP-parametric} is solved repeatedly in a combination of shrinking and receding horizon fashion by varying the two parameters $(h, b_{act})$.

As it can be observed in Algorithm~\ref{alg:mpc}, the proposed MPC algorithm contains two different phases, i.e., a shrinking horizon phase and a receding horizon phase, which are graphically illustrated in Figure~1. The shrinking horizon phase consisting of $l-1$ problems and the receding horizon phase are explained in detail next.

\begin{algorithm}[htb]
\caption{MPC algorithm for switched systems with MTC} 
\label{alg:mpc} 
\begin{algorithmic}[1]
    \Input{The MTC $l$}
    \Initialize{Initialize $t_{act}\gets 0$, $h\gets l$, $b_{act}\gets\emptyset$ }
    \For{$i = 0,1,\ldots$}
	    \State NLP \#1: Solve $P_i(h, b_{act})$ to obtain \revv{$\hat{\mathbf{p}}_i$} given $\tilde{x}_{0|i}$
	    \State CIA: Obtain \revv{$\mathbf{p}_i$} given \revv{$\hat{\mathbf{p}}_i$} using the modified SUR \eqref{SUR}
	    \State NLP \#2: Solve $P_i(h, b_{act})$ to obtain $\mathbf{x}_i,\mathbf{u}_i$ given $\tilde{x}_{0|i}$ and the fixed \revv{$\mathbf{p}_i$}
	    \State $h\gets h-1$
	    \If{$h=0$}
	        \State $h\gets l$
	        \State $b_{act}\gets\emptyset$
	    \Else
	        \State $b_{act}\gets \revv{\mathbf{p}_{0|i}}$
	    \EndIf
	    \State \rev{Return} $u_{0|i}, \mathbf{p}_{0|i}$
    \EndFor
\end{algorithmic}
\end{algorithm}

\begin{figure}[htb]
	\centering
	\includegraphics[width=1\linewidth]{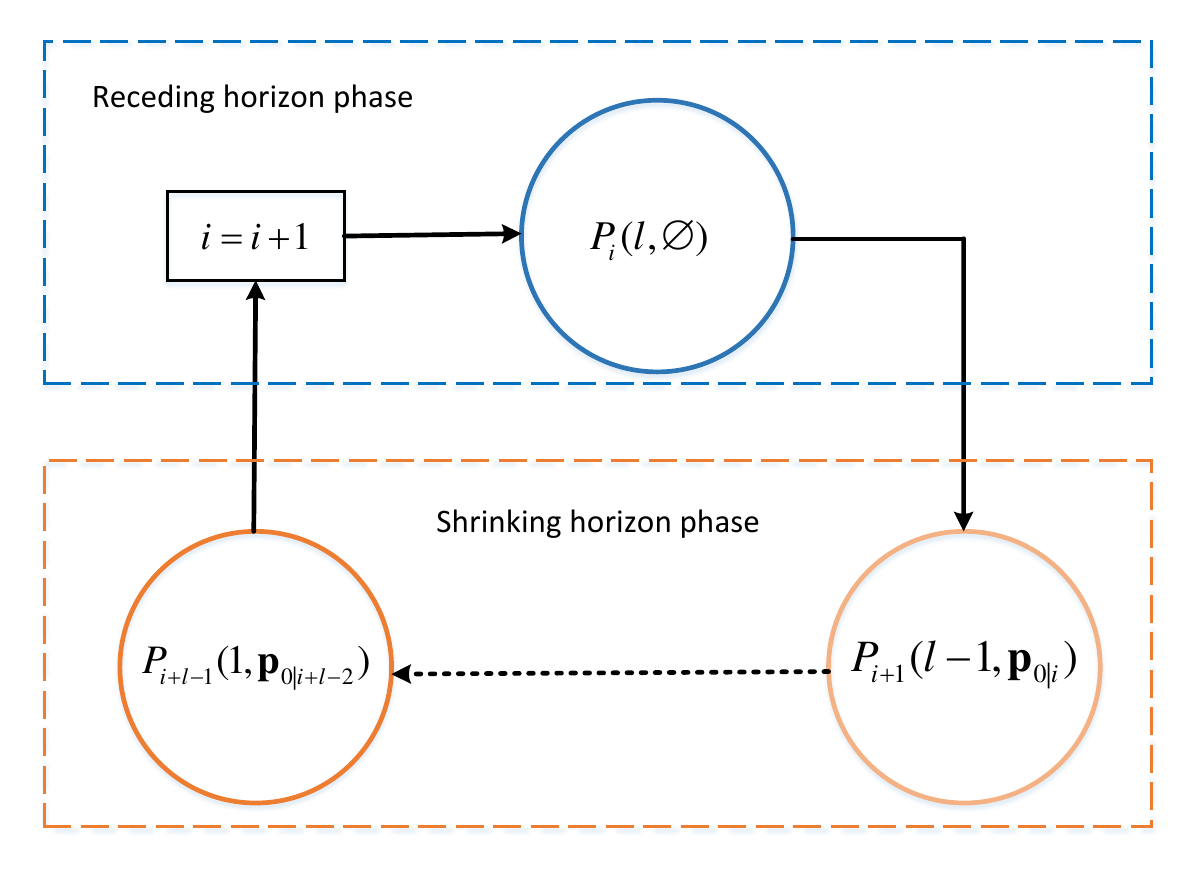}
	\caption{Illustration of shrinking and receding horizon phases for Algorithm~\ref{alg:mpc}.}
	\label{fig:MPC_diagram}
\end{figure}

\subsubsection{Shrinking horizon phase}
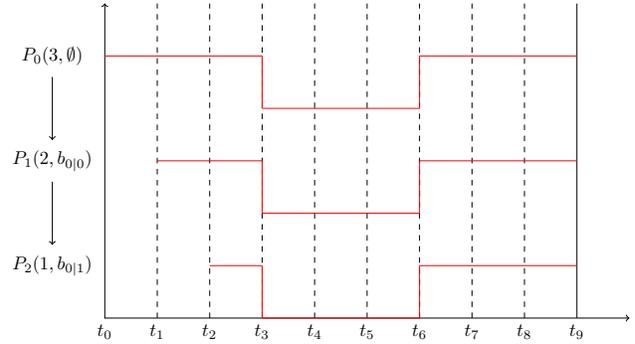
\begin{figure}[htb]
\centering
\begin{adjustbox}{max width=1\linewidth}{
\begin{tikzpicture}
\draw [<->](0,6)--(0,0)--(10,0);
\node [below] at (0,0) {$t_0$};
\draw [dashed,ultra thin] (1,0)--(1,6);
\draw [dashed,ultra thin] (2,0)--(2,6);
\draw [dashed,ultra thin] (3,0)--(3,6);
\node [below] at (1,0) {$t_{1}$};
\node [below] at (2,0) {$t_{2}$};
\node [below] at (3,0) {$t_{3}$};

\draw [dashed,ultra thin] (4,0)--(4,6);
\draw [dashed,ultra thin] (5,0)--(5,6);
\draw [dashed,ultra thin] (6,0)--(6,6);
\node [below] at (4,0) {$t_{4}$};
\node [below] at (5,0) {$t_{5}$};
\node [below] at (6,0) {$t_{6}$};

\draw [dashed,ultra thin] (7,0)--(7,6);
\draw [dashed,ultra thin] (8,0)--(8,6);
\draw [ultra thin] (9,0)--(9,6);
\node [below] at (7,0) {$t_{7}$};
\node [below] at (8,0) {$t_{8}$};
\node [below] at (9,0) {$t_{9}$};

\draw [line width=0.5mm, red] (0,5)--(3,5);
\draw [line width=0.5mm, red] (3,4)--(6,4);
\draw [line width=0.5mm, red] (6,5)--(9,5);

\draw [line width=0.5mm, red] (3,4)--(3,5);
\draw [line width=0.5mm, red] (6,4)--(6,5);

\draw [line width=0.5mm, red] (1,3)--(3,3);
\draw [line width=0.5mm, red] (3,2)--(6,2);
\draw [line width=0.5mm, red] (6,3)--(9,3);

\draw [line width=0.5mm, red] (3,2)--(3,3);
\draw [line width=0.5mm, red] (6,2)--(6,3);

\draw [line width=0.5mm, red] (2,1)--(3,1);
\draw [line width=0.5mm, red] (3,0)--(6,0);
\draw [line width=0.5mm, red] (6,1)--(9,1);

\draw [line width=0.5mm, red] (3,0)--(3,1);
\draw [line width=0.5mm, red] (6,0)--(6,1);

\node [xshift= -1cm] at (0,5) {$P_0(3,\emptyset)$ };
\node [xshift= -1cm] at (0,3) {$P_1(2,\mathbf{p}_{0|0})$ };
\node [xshift= -1cm] at (0,1) {$P_2(1,\mathbf{p}_{0|1})$ };

\draw[->] (-1,4.6)--(-1,3.4);
\draw[->] (-1,2.6)--(-1,1.4);

\end{tikzpicture}
}
\end{adjustbox}
\caption{An illustration of the shrinking horizon phase from one sampling instant to the next. In this example, we start from $i=0$. The number of grid points at the beginning is $N=9$. The two modes are illustrated by flipping the red line. The dwell time constraint is $l=3$.}
\label{fig:shift-block}
\end{figure}

Algorithm~\ref{alg:mpc} starts by solving $P_i(l,\emptyset)$ which consists of $N$ grid points with $M$ blocks of length $l$ \revv{with $h=l$ in \eqref{MPC input block}}. At the next sample, $P_{i+1}(l-1,\mathbf{p}_{0|i})$ is solved where $\mathbf{p}_{0|i}$ is the optimal mode that was computed and implemented at the previous sample. This process is repeated until the sample $i+l-1$ when $P_{i+l-1}(1,\mathbf{p}_{0|i+l-2})$ is solved and the first block has only one interval. As a consequence, a series of problems are solved in the following order:
\begin{align*}\label{problem sequence shrink}
    &P_i(l,\emptyset)\rightarrow P_{i+1}(l-1,\mathbf{p}_{0|i}) \rightarrow  \cdots \rightarrow P_{i+l-1}(1,\mathbf{p}_{0|i+l-2}).
\end{align*}
An illustrative diagram of this shrinking horizon strategy is shown in Fig.~\ref{fig:shift-block}.


\subsubsection{Receding horizon phase}
\begin{figure}[htb]
\centering
\begin{adjustbox}{max width=1\linewidth}{
\begin{tikzpicture}
\draw [<->](0,5)--(0,0)--(12,0);
\node [below] at (0,0) {$t_0$};
\draw [dashed,ultra thin] (1,0)--(1,5);
\draw [dashed,ultra thin] (2,0)--(2,5);
\draw [dashed,ultra thin] (3,0)--(3,5);
\node [below] at (1,0) {$t_{1}$};
\node [below] at (2,0) {$t_{2}$};
\node [below] at (3,0) {$t_{3}$};

\draw [dashed,ultra thin] (4,0)--(4,5);
\draw [dashed,ultra thin] (5,0)--(5,5);
\draw [dashed,ultra thin] (6,0)--(6,5);
\node [below] at (4,0) {$t_{4}$};
\node [below] at (5,0) {$t_{5}$};
\node [below] at (6,0) {$t_{6}$};

\draw [dashed,ultra thin] (7,0)--(7,5);
\draw [dashed,ultra thin] (8,0)--(8,5);
\draw [dashed,ultra thin] (9,0)--(9,5);
\node [below] at (7,0) {$t_{7}$};
\node [below] at (8,0) {$t_{8}$};
\node [below] at (9,0) {$t_{9}$};

\draw [dashed,ultra thin] (10,0)--(10,5);
\draw [dashed,ultra thin] (11,0)--(11,5);
\draw [ultra thin] (12,0)--(12,5);
\node [below] at (10,0) {$t_{10}$};
\node [below] at (11,0) {$t_{11}$};
\node [below] at (12,0) {$t_{12}$};

\draw [line width=0.5mm, red] (0,4)--(3,4);
\draw [line width=0.5mm, red] (3,3)--(6,3);
\draw [line width=0.5mm, red] (6,4)--(9,4);

\draw [line width=0.5mm, red] (3,3)--(3,4);
\draw [line width=0.5mm, red] (6,3)--(6,4);

\draw [line width=0.5mm, red] (3,0)--(6,0);
\draw [line width=0.5mm, red] (6,1)--(9,1);
\draw [line width=0.5mm, red] (9,0)--(12,0);

\draw [line width=0.5mm,  red] (6,0)--(6,1);
\draw [line width=0.5mm, red] (9,0)--(9,1);

\node [xshift= -1cm] at (0,4) {$P_0(3,\emptyset)$ };
\node [xshift= -1cm] at (0,1) {$P_3(3,\emptyset)$ };

\draw[dashed, ->] (-1,3.6)--(-1,1.4);

\end{tikzpicture}
}
\end{adjustbox}
\caption{An illustration of the receding horizon phase using the same settings as in Fig. \ref{fig:shift-block}. The problems between $P_0(3,\emptyset)$ and $P_3(3,\emptyset)$ are $P_1(2,\mathbf{p}_{0|0})$ and $P_2(1,\mathbf{p}_{0|1})$ in Fig. \ref{fig:shift-block}.}
\label{fig:rhc-block}
\end{figure}
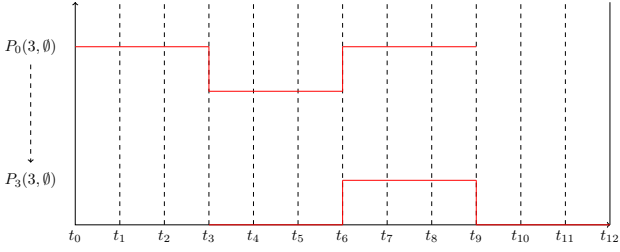

The receding horizon phase is performed after the shrinking horizon phase has ended, which is triggered when $h=0$. In this phase, problem $P_{i+l}(l,\emptyset)$ is formulated and solved by introducing an additional block of length $l$ at the tail of the prediction horizon to recover the original prediction length. This is equivalent to shifting the problem $P_i(l,\emptyset)$ $l$ steps forward. An illustrative diagram of this receding horizon step is shown in Fig.~\ref{fig:rhc-block}. 

\section{Recursive feasibility guarantees}
In this section, we first show that Algorithm~\ref{alg:mpc} is recursively \revv{$\delta$-feasible} by using a $l$-step control invariant ($l$-CI) terminal set $\mathcal{X}_f$. Then we show how to compute such sets in a tractable way for switched linear systems. We focus on recursive feasibility guarantees because this is the most important property that an MPC algorithm must satisfy, to facilitate implementation in practice. Future work, which is beyond the scope and page limits of the current paper, will also consider developing closed-loop stability guarantees for the MPC controller generated by Algorithm~\ref{alg:mpc}, by means of non-monotonic Lyapunov functions \cite{lazar2013stability}.

\begin{defi}\label{def: l-step invariance}
	The set $\mathcal{X}_f$ is an admissible $l$-step control invariant ($l$-CI) for the dynamics $\phi(\cdot,\cdot,\cdot)$ if for all $x_0\in\mathcal{X}_f$, there exists a control sequence $\{u_{0}, u_1, \ldots,u_{l-1}\}\in\mathcal{U}^{l}$ and a mode mapping $\mathbf{p}_0\in\{0,1\}^Q,\, \sum_{j=1}^Q p_{0}^{j}=1$ consisting of $l$ consecutive blocked mode sequences such that
	\begin{equation}\label{nonlinear discrete dynamics}
	\begin{aligned}
	&x_{i+1}=\phi(x_{i},u_{i},\mathbf{p}_0),\quad\forall i=0,\ldots,l-1,\\
	&x_i\in\mathcal{X},\quad\forall i=0,\ldots,l-1,\\
	& \revv{g(x_i,u_i)\leq 0, \quad\forall i=0,\ldots,l-1,}\\
	&x_l\in \mathcal{X}_f.
	\end{aligned}
	\end{equation}
\end{defi}
Next, we define a mapping $\Phi(\cdot,\cdot,\cdot)$ which maps the initial state $x_0$ to an $l$-step ahead state $x_l$ via an admissible continuous control input sequence $\{u_{0}, u_1, \ldots,u_{l-1}\}$ and mode mapping $\mathbf{p}_0$, i.e.
\[x_l=\Phi(x_0, \{u_{0}, u_1, \ldots,u_{l-1}\}, \mathbf{p}_0),\quad \forall l\in\mathbb{N}.\]

\begin{assmp}\label{assmp:terminal set is l-CI}
	The terminal set $\mathcal{X}_f$ in \eqref{NLP-parametric} is an admissible $l$-CI set for the dynamics $\phi$, with $l$ equal to the dwell time interval length.
\end{assmp}

Assumption \ref{assmp:terminal set is l-CI} is inspired by the invariant $(k,\lambda)$ contractive set which defines set invariance in $k$ steps \cite{lazar2013stability}. Assumption \ref{assmp:terminal set is l-CI} also adopts a similar idea employed in event-triggered MPC \cite{varutti2009event} where the set invariance is defined at a future time point.


The following theorem states the main result of this section.
%
\begin{thm}\label{thm2}
If Assumption~\ref{assmp:feasible dwell time} and Assumption~\ref{assmp:terminal set is l-CI} hold, the MPC scheme presented in Algorithm~\ref{alg:mpc} is recursively \rev{$\delta$-feasible}, i.e. if $P_0(l,\emptyset)$ is feasible, \rev{problems $P_{i}(l,b_{act})$ for all $i\in\mathbb{N}$ and corresponding admissible pairs $(l,b_{act})$ are $\delta$-feasible.}
\end{thm}
\begin{pf}
The proof consists of two parts: \emph{i)} proving \revv{$\delta$-feasibility} of the NLP \#2 given a feasible solution obtained from the NLP \#1 and the SUR \eqref{SUR}, at the current sampling instant with the same initial condition; \emph{ii)} given a \revv{$\delta$-feasible} NLP \#2 at the current sampling instant, proving feasibility of the NLP \#1 at the next sampling instant 

\emph{i)}: We start by assuming the NLP \#1 of problem $P_i(l,\emptyset)$ is feasible, for an arbitrary $i\in\mathbb{N}$. Based on Assumption~\ref{assmp:feasible dwell time},  the NLP \#2 of problem $P_i(l,\emptyset)$ is \rev{$\delta$-feasible}, with the same initial condition.

\emph{ii)}: Next, we prove the NLP \#1 at the next sampling instant is feasible. 

(A) The shrinking horizon phase: 

The optimal solutions to the NLP \#2 of problem $P_i(l,\emptyset)$ are defined by $\{u_{0|i}, u_{1|i}, \ldots,u_{N-1|i}\}\in\mathcal{U}^{N}$ and \revv{$\{\mathbf{p}_{0|i}, \mathbf{p}_{1|i}, \ldots, \mathbf{p}_{M-1|i}\}\in\{0,1\}^{Q\times M}$} satisfying \eqref{MPC input block}. The solutions lead to a $\delta$-feasible state trajectory $\{x_{0|i}, x_{1|i}, \ldots,x_{N|i}\}$ \revv{according to Assumption~\ref{assmp:feasible dwell time}}. The optimal mode that is fed back to the system is \revv{$b_{act}:=\mathbf{p}_{0|i}$}. At sample $i+1$, the prediction horizon is shrunk and the NLP \#1 of problem \revv{$P_{i+1}(l-1,\mathbf{p}_{0|i})$} is solved. The solutions denoted by
\begin{subequations}
\begin{align}
    &\{u_{0|i+1}, \ldots,u_{N-2|i+1}\}=\{u_{1|i}, \ldots,u_{N-1|i}\},\\
    &\revv{\mathbf{p}_{0|i+1}=\mathbf{p}_{0|i}=b_{act},}\\
    &\revv{\{\mathbf{p}_{1|i+1},\ldots,\mathbf{p}_{M-1|i+1}\}=\{\mathbf{p}_{1|i},\ldots,\mathbf{p}_{M-1|i}\},}
\end{align}
\end{subequations}
lead to a state trajectory 
\begin{equation}
    \{x_{0|i+1},\ldots,x_{N-1|i+1}\}=\{x_{1|i},\ldots,x_{N|i}\},
\end{equation}
which is feasible in terms of dynamics $\phi(\cdot,\cdot,\cdot)$ and constraints specified by $g(\cdot,\cdot)$. \revv{Note that from NLP \#2 to NLP \#1 there is no integer approximation error hence it is not necessary to use the notion of $\delta$-feasibility. Here, both $\mathbf{p}_{0|i+1}$ and $\mathbf{p}_{0|i}$ are using the value of $b_{act}$, but they incorporate a different number of blocks due to horizon shrinking, i.e.
\begin{align*}
    &\mathbf{p}_{0|i+1} \gets \mathbf{b}_1=\ldots=\mathbf{b}_{h-1},\\
    &\mathbf{p}_{0|i} \gets \mathbf{b}_0=\ldots=\mathbf{b}_{h-1}.
\end{align*}}
By combining \emph{i)} and \emph{ii)}, recursive $\delta$-feasibility can be constructed until sample $i+l-1$ when problem $P_{i+l-1}(1,\mathbf{p}_{0|i+l-2})$ is solved, for every $i$ starting from $0$ with an increment $l$.

(B): The receding horizon phase:

Given that the NLP \#2 of problem $P_{i+l-1}(1,\mathbf{p}_{0|i+l-2})$ is $\delta$-feasible, we prove the NLP \#1 of problem $P_{i+l}(l,\emptyset)$ is feasible. Choose the solution of the first $N-l$ intervals and \revv{$M-1$ blocks} as
\begin{subequations}
\begin{align}
    &\{u_{0|i+l}, \ldots,u_{N-l-1|i+l}\}=\{u_{1|i+l-1}, \ldots,u_{N-l|i+l-1}\},\\
    & \revv{\{\mathbf{p}_{0|i+l},\ldots,\mathbf{p}_{M-2|i+l}\}=\{\mathbf{p}_{1|i+l-1},\ldots,\mathbf{p}_{M-1|i+l-1}\}}
\end{align}
\end{subequations}
which leads to a state trajectory
\begin{subequations}
\begin{align}
    &\{x_{0|i+l},\ldots,x_{N-l-1|i+l}\}=\\
    & \quad \quad \{x_{1|i+l-1},\ldots,x_{N-l|i+l-1}\}\in\mathcal{X}^{N-l},\nonumber\\
    &x_{N-l|i+l} = x_{N-l+1|i+l-1}\in\mathcal{X}_f
\end{align}
\end{subequations}
and such that 
\begin{equation}
    g(x_{k|i+l},u_{k|i+l})\leq 0,\,k=0,1,\ldots,N-l-1.
\end{equation}
As a result, the appended block at the the tail of the prediction horizon starts with the initial state $x_{N-l|i+l}\in\mathcal{X}_f$. \revv{According to Assumption \ref{assmp:terminal set is l-CI}, there exist a binary block representation 
\begin{equation*}
    \mathbf{p}_{M-1|i+l}\gets \mathbf{b}_{h+(M-2)l|i+l}=\ldots=\mathbf{b}_{h+(M-1)l|i+l},
\end{equation*}
and a sequence of continuous inputs $$\{u_{N-l|i+l}, \ldots,u_{N-1|i+l}\}\in\mathcal{U}^{l},$$} such that 
\begin{subequations}
\begin{align}
    x_{N-l+k|i+l}&=\phi(x_{N-l|i+l},u_{N-l|i+l},p_{M-1|i+l})\in\mathcal{X},\nonumber\\
    &\forall k=1,\ldots,l-1,\\
    x_{N|i+l}&=\Phi(x_{N-l|i+l}, \{u_{N-l|i+l}, \ldots,u_{N-1|i+l}\},\nonumber\\
    &p_{M-1|i+l})\in\mathcal{X}_f,
\end{align}
\end{subequations}
and 
\begin{equation}
    g(x_{N-l+k|i+l},u_{N-l+k|i+l})\leq 0,\,\forall k=0,\ldots,l-1.
\end{equation}
Therefore, the NLP \#1 of problem $P_{i+l}(l,\emptyset)$ admits a feasible solution. Since $i\in\mathbb{N}$ was arbitrarily chosen, this completes the proof. \hfill{$\qed$}
\end{pf}

\begin{rmk}
Theorem~\ref{thm2} requires Assumption~\ref{assmp:feasible dwell time} and \ref{assmp:terminal set is l-CI} which may not hold in practice. In particular, Assumption~\ref{assmp:feasible dwell time} requires the MTC to be small enough because $\delta_1,\delta_3$ in \eqref{defi1 equa} cannot be made arbitrary small when MTC is a hard constraint for the OCP \eqref{NLP-blocked}. Hence, Assumption~\ref{assmp:feasible dwell time} is reasonable by considering that the controlled system losses the degree of freedom to switch between modes if the MTC is too long. Assumption~\ref{assmp:terminal set is l-CI} is more restrictive in the sense that it imposes constraints on the terminal $l$-CI set $\mathcal{X}_f$, \revv{and requires constraint fulfillment of the general nonlinear constraint $g(\cdot,\cdot)$, which is hard to verify a priori.} Note however that Assumption~\ref{assmp:terminal set is l-CI} employs a relaxed version of the usual control invariant terminal set, which is only required to be periodically control invariant. This is less conservative compared to the standard control invariant terminal set condition typically used in nonlinear MPC to establish recursive feasibility.
\end{rmk}


\subsection{Computation of the SRCI terminal set}

The explicit computation of $\mathcal{X}_f$ that satisfies Assumption \ref{assmp:terminal set is l-CI} is not straightforward, even for standard nonlinear MPC algorithms \cite{mayne2000constrained}. In this work, we develop an iterative algorithm to compute a specific type of switch-robust $l$-CI set for linear switched systems. To this end consider the following definition. 
\begin{defi}\label{def:l-SRCI}
	The set $\mathcal{X}_f$ is an admissible $l$-step switch-robust CI (l-SRCI) set for the dynamics $\phi(\cdot,\cdot,\cdot)$ if for all $x_0\in\mathcal{X}_f$, there exists $\{u_{0}, u_1, \ldots,u_{l-1}\}\in\mathcal{U}^{l}$ such that for all $p_0\in\{0,1\}^Q$, $\sum_{j=1}^Q p_{0}^{j}=1$, it holds that
	\begin{equation}
	\begin{aligned}
	&x_{l}=\phi(x_{l-1},u_{l-1},\mathbf{p}_0)\\
	 & \quad=\Phi(x_0, \{u_{0}, u_1, \ldots,u_{l-1}\}, \mathbf{p}_0) \in\mathcal{X}_f,\\
	&\revv{x_i\in\mathcal{X},\quad\forall i=0,\ldots,l-1,}\\
	& \revv{g(x_i,u_i)\leq 0, \quad\forall i=0,\ldots,l-1,}\\
	\end{aligned}
	\end{equation}
\end{defi}

\begin{prop}\label{thm:l-SRCI to l-CI}
If $\mathcal{X}_f$ is an admissible $l$-SRCI set, then it is an admissible $l$-CI set.
\end{prop}
The claim of Proposition \ref{thm:l-SRCI to l-CI} directly follows based on Definition \ref{def: l-step invariance} and \ref{def:l-SRCI}.  Therefore, if we can compute a terminal set $\mathcal{X}_f$ as an admissible $l$-SRCI set, this is a sufficient condition for Assumption \ref{assmp:terminal set is l-CI} to hold. 

\revv{In what follows we consider the case when the general nonlinear constraint $g(\cdot,\cdot)$ is not present. This is justified by the common occurrence of polytopic state and input constraints in practical applications, case in which set computations can be implemented.} Algorithm~\ref{alg:srci} summarizes the computation of a $l$-SRCI set in this case. Therein, for all $j=1,\ldots,Q$, the predecessor set is defined by
\begin{align}
    \text{Pre}(\mathcal{X};j)=&\{x\in\mathbb{R}^{n_x}:\,\exists u\in\mathcal{U},\,\phi(x,u,e^j)\in \mathcal{X}\},\\
    \text{Pre}^{k+1}(\mathcal{X};j)=&\{x\in\mathbb{R}^{n_x}:\,\nonumber\\&\exists u\in\mathcal{U},\phi(x,u,e^j)\in \text{Pre}^k(\mathcal{X};j)\},
\end{align}
where $\text{Pre}^0(\mathcal{X};j)=\mathcal{X}$ and $e^j$ is defined by
\begin{equation}\label{mode p_j}
    e^j=[0,\ldots,0,1,0,\ldots,0]^\top\in\{0,1\}^Q,
\end{equation}
with all elements zero except the $j$-th element.

\begin{algorithm}[htb]
\caption{$l$-SRCI set computation \revv{without constraints $g(\cdot,\cdot)$}} 
\label{alg:srci} 
\begin{algorithmic}[1]
    \Input{$\mathcal{X},l$}
    \Initialize{Initialize $\mathcal{X}_f^0=\mathcal{X}, i=0$ }
    \Repeat
	    \State Update $\mathcal{X}_f^{i+1}=\mathcal{X}_f^{i}\cap_{j\in\{1,\ldots,Q\}} \text{Pre}^l(\mathcal{X}_f^i;j)$
	    \State $i=i+1$
	\Until{$\mathcal{X}_f^{i+1}=\mathcal{X}_f^{i}$}
	\State $\mathcal{X}_f=\mathcal{X}_f^{i}$
\end{algorithmic}
\end{algorithm}

Algorithm~\ref{alg:srci} updates $\mathcal{X}_f$ by taking the intersection of the $l$-step predecessor sets for all modes. Proving convergence of Algorithm~\ref{alg:srci} in general (i.e., for nonlinear dynamics) is a non-trivial, open problem. However, consider the case of a switched linear autonomous system, i.e., 
\begin{equation}
\label{sys:switched}
    x_{k+1} = \sum_{j=1}^Q \left(b_{k,j}A_j\right) x_k. 
\end{equation}
The $l$-step terminal state can be computed as
\begin{align}\label{computation of x_l}
    x_{l} &= \left( \sum_{j=1}^Q p_{0}^{j}\: A_j \right)^l x_0,
\end{align}
where $p_0$ represents the blocked binary variables over $l$ steps. In this case, the $l$-step predecessor set can be considered as a $1$-step set using
\begin{align}\label{l-step pre set}
    \text{Pre}^l(\mathcal{X};j)=&\{x\in\mathbb{R}^{n_x}:\,\tilde{A}_jx\in \mathcal{X}\},
\end{align}
where $\tilde{A}_j=(A_j)^l$. 

\begin{rmk}
For system with dynamics \eqref{computation of x_l} and assuming $\mathcal{X}_f$ compact with $0\in\text{int}(\mathcal{X}_f)$, Algorithm~\ref{alg:srci} terminates in finite number of steps if $\tilde{A}_j,\,j=1,\ldots,Q$ are stable \cite{lazar2006stabilizing}. If at least one of $\tilde{A}_j$ is stable, it is likely that Algorithm~\ref{alg:srci} asymptotically converges to the set $\{0\}$ due to intersections of predecessor sets. 
\end{rmk}

\begin{rmk}
Consider linear systems with dynamics
\begin{align}\label{linear dynamics}
    x_{k+1} = \sum_{j=1}^Q b_{k}^{j}\left(A_j x_k + B_j u_k\right). 
\end{align}
If there exits a feedback law $u_k = K_jx_k$ such that $A_j+B_jK_j,\,j=1,\ldots,Q$ are stable, then we can apply Algorithm~\ref{alg:srci} to \eqref{linear dynamics} using the predecessor set defined as \eqref{l-step pre set} where $\tilde{A}_j=(A_j+B_jK_j)^l$ \cite{lazar2006stabilizing}.
\end{rmk}

\begin{rmk}
Definition \ref{def:l-SRCI} is similar to the SRCI set defined in \cite{danielson2019necessary} in the way that the latter defines different feasible sets under different modes.  The SRCI set notion defined in \cite{danielson2019necessary} is different and requires one step invariance for one mode and $l$-step invariance for others. A similar algorithm was presented in \cite{danielson2019necessary} for computing corresponding $1$-step SRCI sets. The adopted $l$-step SRCI definition in this paper requires $l$-step invariance under all modes and reduces to a subset of the SRCI definition in \cite{danielson2019necessary} only when $l=1$. 
\end{rmk}

\revv{
\begin{rmk}
In case that the general nonlinear constraints specified by $g(\cdot,\cdot)$ are present, Algorithm~\ref{alg:srci} cannot guarantee the fulfillment of such constraints explicitly, neither provides a feasible way to construct the corresponding $l$-SRCI set. Yet, the constraints $g(\cdot,\cdot)$ could be enforced implicitly by computing a polytopic inner approximation of $\mathcal{G}:=\{(x,u)\in\mathcal{X}\times\mathcal{U}\ : \ g(x,u)\leq 0\}$ and constraining (or scaling) $\mathcal{X}_f$ such that $\mathcal{X}_f\subseteq\text{Proj}_\mathcal{X}(\mathcal{G})$ ($\text{Proj}$ denotes a projection operator). It is worth to mention also that computing control invariant sets in the presence of nonlinear state and input constraints is a non--trivial open problem even for standard nonlinear MPC, i.e., without MTC present. 
\end{rmk}
}

\begin{table*}[!htb]
\centering
\caption{\revv{Accumulated objective function value,} constraint violation and computation time comparison for Example 1. A1 stands for Algorithm~\ref{alg:mpc}. CPT denotes the average computational time in milliseconds per sampling instant.}
\label{tab:ex1_comparison}
\rev{
\begin{tabular}{ccc|cc|ccc}
\hline
    & \multicolumn{2}{c|}{$E$} & \multicolumn{2}{c|}{$res$} & \multicolumn{3}{c}{CPT~[ms]}         \\ \cline{2-8} 
    & A1      & (MPC+CIA)  & A1       & (MPC+CIA)    & A1    & (MPC+CIA)$_1$ & (MPC+CIA)$_2$ \\
$l=2$ & \revv{6.433}  &  \revv{6.545}  & 0.267    & 0.064         & 0.722 & 23.34         & 1.93      \\
$l=4$ &  \revv{6.566}  &  \revv{6.566}      & 0.199    & 0.199         & 0.581 & 22.84         & 1.88      \\
$l=5$ &  \revv{6.001}  &   \revv{6.543}     & 0.411    & 0.411         & 0.540 & 21.80         & 1.78      \\ \hline
\end{tabular}
}
\end{table*}

\section{Numerical Examples}\label{sec5}
In this section, two numerical examples are presented to show the effectiveness and efficiency of Algorithm~\ref{alg:mpc}. Comparisons have been made with \rev{two variants} of the algorithm proposed in \cite{burger2019design} which employs the decomposition method described in Section 2.2. \rev{We denote this algorithm as (MPC+\rev{CIA})$_1$ and (MPC+\rev{CIA})$_2$, which uses an MILP and a tailored branch-and-bound (BnB) algorithm \cite{burger2019design} for the CIA problem, respectively. The tailored BnB algorithm has been reported to reduce the computation time for the CIA problem significantly.} The simulations were performed on a \rev{PC running Windows 10 with Intel i5 8500 CPU at 3.0GHz.} The NLP in the form of \eqref{NLP-blocked} is solved in \software{MATMPC} \cite{chen2019matmpc} which is an open-source \software{MATLAB} based nonlinear MPC software. The MILP is solved using \software{CPLEX} in \software{MATLAB} \rev{and the tailored BnB implementation is taken from the toolbox \software{pycombina} \cite{burger2019design}.}

\subsection{Example 1}
\begin{figure}[htb]
    \centering
    \includegraphics[width=.9\linewidth]{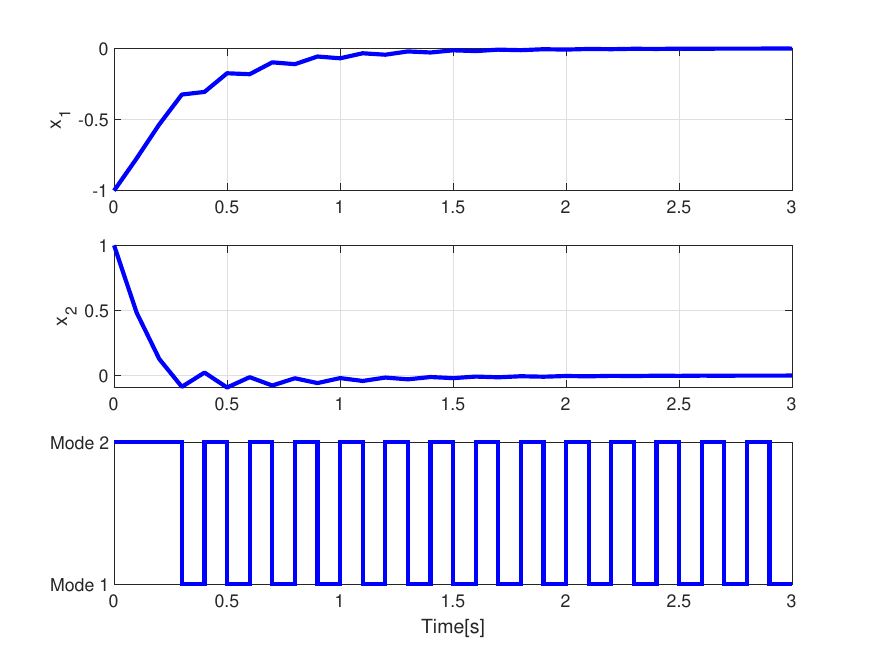}
    \caption{State and mode trajectories of problem \eqref{ex1:dynamics} without MTC.}
    \label{fig:ex1_noMTC}
\end{figure}

Consider a linear autonomous system governed by the following dynamics
\begin{align}\label{ex1:dynamics}
    \begin{cases}
    \dot{x}(t) = \begin{bmatrix}
    -5 & -3 \\
    5 & -1
    \end{bmatrix}x(t)\quad &\text{Mode 1},\\
    \dot{x}(t) = \begin{bmatrix}
    -1 & 5 \\
    -3 & -5
    \end{bmatrix}x(t)\quad &\text{Mode 2}.
    \end{cases}
\end{align}
We design an MPC controller to regulate the states from the initial state $x(0)=(-1,1)^\top$ to the origin $x_T=(0,0)$. The MPC controllers used in this simulation are configured as follows. The cost function for the on-line optimization problem is given by
\begin{align}\label{ex1:cost}
    J(x)=\sum_{k=0}^{N-1}x_k^\top Q x_k + x_N^\top P x_N
\end{align}
where $Q=[1~0;0~1], P=[10~0;0~10]$. The number of grid points is $N=20$ and the sampling time is $\Delta t=0.1$~s. The states are constrained in the set $\mathcal{X}=\{(x_1,x_2):-1\leq x_1\leq 0.05, -0.05\leq x_2\leq 1\}$. 
Fig.~\ref{fig:ex1_noMTC} shows the state and mode trajectories obtained using Algorithm~\ref{alg:mpc} with $l=1$, i.e. without MTC. The controller needs to switch between the two modes frequently to regulate both states to the origin to avoid violating the state constraints. We define two criteria to evaluate the qualitative performance of the algorithms as 
\begin{equation}
    \begin{aligned}
        \revv{E =} & \revv{\sum_{i=0}^I J(x(t_i)),}\\
        res = & \sum_{i=0}^{\revv{I}}\left\{\max(x_1(t_i)-0.05,0)+\max(x_2(t_i)-1,0)+\right.\\
        &\left.\max(-x_1(t_i)-1,0)+\max(-x_2(t_i)-0.05,0)\right\},
    \end{aligned}
\end{equation}
where $E$ is the sum of \revv{objective function values} over the entire closed-loop simulation. $res$ is the accumulated violation of state constraints, and \revv{$i\in\mathbb{N}$} is the sampling instant, \revv{whose range $I$ depends on the duration of the closed-loop simulation}. The results of the three tested algorithms are shown in Table~\ref{tab:ex1_comparison} for different MTC $l$. The two (MPC+CIA) variants \revv{employ} the same MPC algorithm with different CIA algorithms, hence their names are not specified when comparing the control performance and constraint violation criteria. It can be observed that Algorithm~\ref{alg:mpc} has \revv{similar or even smaller $E$ value, but has higher constraint violation} when $l=2$. However, Algorithm~\ref{alg:mpc} is much faster than (MPC+CIA) in all cases, largely due to the fact that the move blocking MPC \eqref{NLP-blocked} is easier to solve and that the SUR step is much cheaper than implementing the MILP or the tailored BnB algorithm. In this example, the constraint is violated because the system is autonomous and the MPC controller is not able to switch that frequently without the help of external control inputs to maintain the states inside the feasible region. Note that for this example there is no need to formulate and solve NLP \#2 using Algorithm~\ref{alg:mpc}.

\begin{figure}[htb]
    \centering
    \includegraphics[width=.9\linewidth]{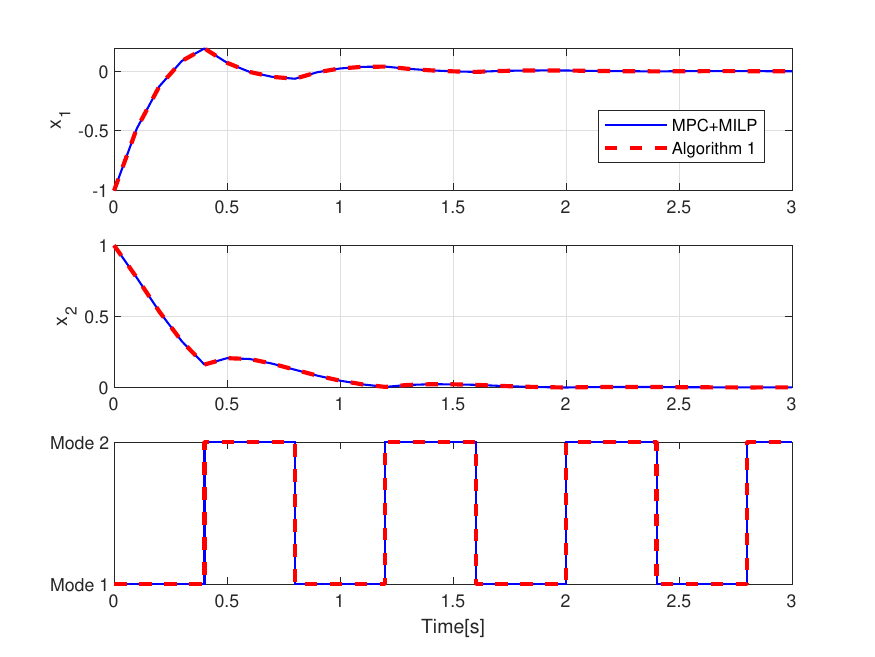}
    \caption{State and mode trajectories of problem \eqref{ex1:dynamics} using the MPC+MILP algorithm and Algorithm \ref{alg:mpc}, with MTC $l=4$.}
    \label{fig:ex1_compare}
\end{figure}

\begin{figure}[htb]
    \centering
    \includegraphics[width=.9\linewidth]{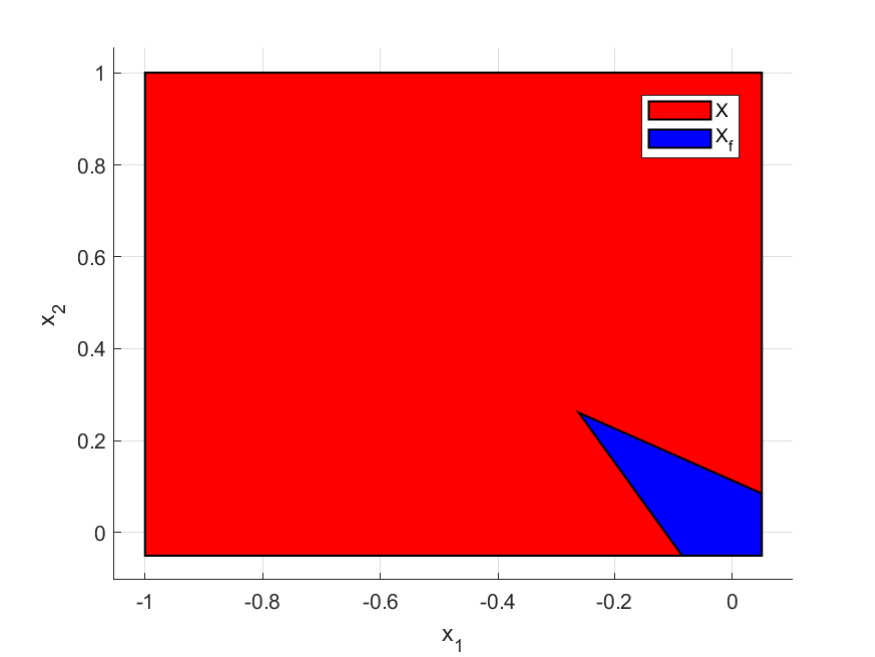}
    \caption{The feasible set $\mathcal{X}$ and the SRCI terminal set $\mathcal{X}_f$ of problem \eqref{ex1:dynamics} when MTC $l=4$, computed by Algorithm 2.}
    \label{fig:ex1_SRCI}
\end{figure}

Fig.~\ref{fig:ex1_compare} shows the state and mode trajectories using Algorithm \ref{alg:mpc} and the (MPC+CIA) algorithm when $l=4$, which requires $\underline{\Delta \tau}=0.4$~s. We calculate the $l$-step SRCI terminal set $\mathcal{X}_f$ using Algorithm 2, which is shown in Fig~. \ref{fig:ex1_SRCI}. 
Consistently with the results reported in Table~\ref{tab:ex1_comparison}, Algorithm~\ref{alg:mpc} and the (MPC+CIA) algorithm result in the same closed-loop trajectories/performance. 

\begin{table*}[!tb]
\centering
\caption{\revv{Scaled objective function value} and computation time comparison for Example 2. A1 stands for Algorithm~\ref{alg:mpc}. MPC stands for the standard MPC using A1 without MTC. CPT denotes the average computational time in millisecond per sampling instant.}
\label{tab:ex2_comparison}
\rev{
\begin{tabular}{cccc|cccc}
\hline
    & \multicolumn{2}{c}{$E$} &     & \multicolumn{3}{c}{CPT~[ms]}         &     \\ \cline{2-8} 
    & A1      & (MPC+CIA)  & MPC & A1    & (MPC+CIA)$_1$ & (MPC+CIA)$_2$ & MPC \\
$l=4$ & \revv{7.034}   &  \revv{2.902}     &  \revv{2.758}   & 30.53 & 73.03         &  43.43     &  50.23   \\
$l=5$ &  \revv{13.441}  & \revv{3.202}     &  -   & 32.63 & 80.81         & 45.80      & -   \\
$l=8$ &  \revv{59.444}  & \revv{4.202}      &   -  & 34.97 & 86.74         & 51.50      &  -   \\ \hline
\end{tabular}
}
\end{table*}

\subsection{Example 2}

Consider a bevel-tip flexible needle system \cite{vasudevan2013consistentp2} governed by the dynamics 

\begin{align}\label{ex2:dynamics}
\begin{cases}
\dot{x}(t)=\begin{bmatrix}\sin x_5(t)u_1(t)\\-\cos x_5(t)\sin x_4(t)u_1(t)\\\cos x_4(t)\cos x_5(t)u_1(t)\\\kappa\cos x_6(t)\sec x_5(t)u_1(t)\\\kappa\sin x_6(t)u_1(t)\\-\kappa\cos x_6(t)\tan x_5(t)u_1(t)\end{bmatrix}\quad &\text{Mode 1},\\ \dot{x}(t)=\begin{bmatrix}0\\0\\0\\0\\0\\u_2(t)\end{bmatrix}\quad &\text{Mode 2}.
\end{cases}
\end{align}
where $\kappa=0.22$ is the curvature of the needle. The first three states are positions of the needle and the last three are the yaw, pitch and roll angle. The needle is pushed in the first mode and is turning in the second. The control input are the insertion speed $u_1$ and the rotation speed $u_2$. The objective of the controller is to drive the needle from $x(0)=(0,0,0,0,0,0)^\top$ to the destination $x_T=(-2,3.5,10,0,0,0)$ while minimizing the energy cost, defined by the function 
\begin{equation}\label{ex2:cost}
J(x,u)=\sum_{k=0}^{N-1}u_{k}^\top Q u_{k} + \rev{10}(x_N-x_T)^\top(x_N-x_T)
\end{equation}
where $Q=0.01I_2$ is the scaled identity matrix. In addition, the needle must avoid three obstacles defined by spheres centered at $(0,0,5),\,(1,3,7),\,(-2,0,10)$ with radius $2$. The control inputs are constrained as $u_1\in [0,5], u_2\in [-\pi/2,\pi/2]$. The performance evaluation is defined as the summation of \revv{the objective function} defined by
\revv{
\begin{equation}
    E =  \sum_{i=0}^I J(x(t_i),u(t_i)),
\end{equation}
where $I$ denotes the range of the closed-loop simulation as in the first example.} For this example, we choose the number of grid points $N=40$ and the sampling time $\Delta t=0.1$ s. Table~\ref{tab:ex2_comparison} shows the performance and computational time of the three algorithms for this example with varying MTC lengths. The average computation time per sampling instant of the Algorithm~\ref{alg:mpc} is much smaller than the (MPC+CIA)$_1$ and is also smaller with less margin than the (MPC+CIA)$_1$ variant. The reason is twofold: i) Algorithm~\ref{alg:mpc} formulates NLP problems with less decision variables and allows tailored move blocking algorithms for solving the NLP \#1. The computational time for this step using Algorithm~\ref{alg:mpc} is around 30~ms while that for (MPC+CIA) is around 45~ms; ii) the SUR step (around 0.2~ms) is computationally cheaper than employing MILP (around 30~ms) or the tailored BnB solver (around 2~ms). A second observation is that the performance of Algorithm~\ref{alg:mpc} is poorer than (MPC+CIA), and their performance gap grows when the length of MTC increases. This is due to the fact that Algorithm~\ref{alg:mpc} losses the degree of freedom to drive the system towards its objective when MTC is not sufficiently small, while (MPC+CIA) maintains its degree of freedom when solving the NLPs and impose MTC only at the CIA step. Both methods are faster but have poorer performance than using MPC without MTC. Note that the computation speedup for this example is not as significant as for Example 1 since the NLPs defined by \eqref{ex2:dynamics} and \eqref{ex2:cost} are much harder to solve than the corresponding MILPs, due to the complex nature of the nonlinear constraints.

\begin{figure}[htb]
     \centering
     \includegraphics[width=\linewidth]{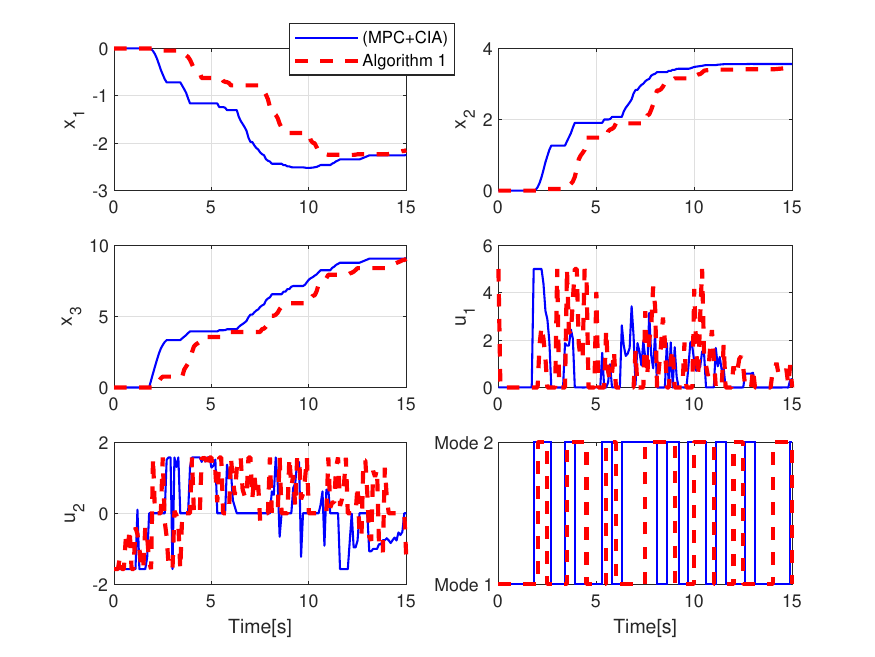}
     \caption{State, control and mode trajectories of problem \eqref{ex2:dynamics}. Algorithm~\ref{alg:mpc} and the (MPC+CIA) algorithm are constrained with MTC $l=5$.}
     \label{fig:ex2_compare}
 \end{figure}

\begin{figure}[htb]
	\centering
	\includegraphics[width=\linewidth]{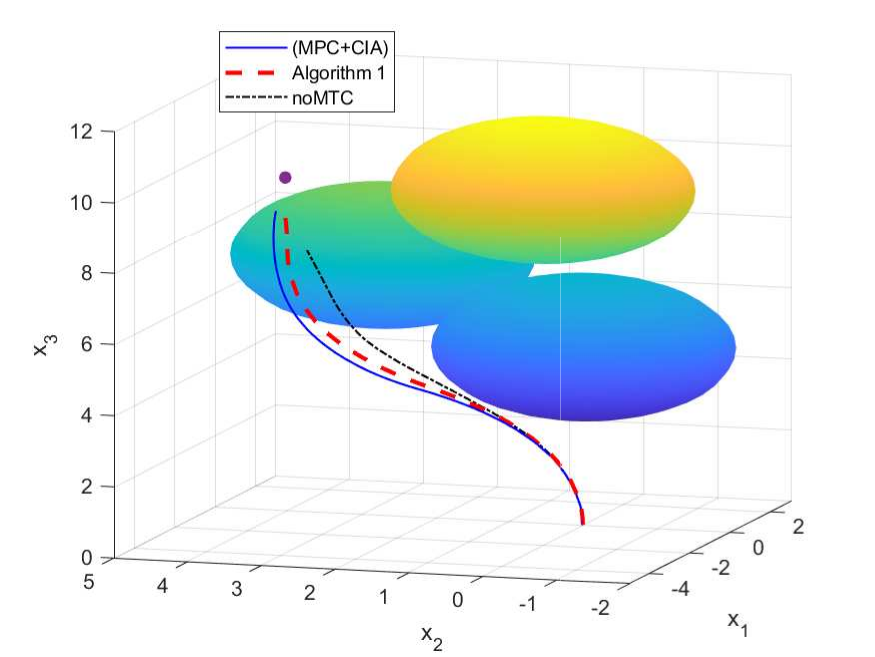}
	\caption{\rev{Closed-loop position trajectories of problem \eqref{ex2:dynamics} using the \revv{two} algorithms. The three obstacles are represented as spheres. The objective point $(-2, 3.5, 10)$ is a filled black circle. Algorithm~\ref{alg:mpc} and the (MPC+CIA) are constrained with MTC $l=5$.}}
	\label{fig:ex2_traj}
\end{figure}

The closed-loop position trajectories when $l=5$ are shown in Fig.~\ref{fig:ex2_traj}. Both Algorithm~\ref{alg:mpc} and the (MPC+CIA) can drive the needle towards the objective point. Fig.~\ref{fig:ex2_compare} shows the state, control and mode trajectories using the \revv{two} algorithms. The state and the mode trajectories show that Algorithm~\ref{alg:mpc} reacts slower in driving the needle to its target compared to (MPC+CIA), because of lower degree of freedom in terms of switching. The control trajectories show that, to achieve the control goal, Algorithm~\ref{alg:mpc} uses more control input power and changes the control input more frequently. When $l=8$, Algorithm~\ref{alg:mpc} fails to drive the needle to the objective point while (MPC+CIA) succeeds.
Hence, Algorithm~\ref{alg:mpc} provides a trade-off between computational cost and control performance, which is consistent with the theoretical results from Assumption~\ref{assmp:feasible dwell time} and Theorem~\ref{thm2}.

\section{Conclusion}\label{sec6}

This paper developed an efficient MPC algorithm for switched systems subject to MTC. Using the decomposition method for MIOC problems, the MTC have been embedded into the first NLP using move blocking. As a result, a simple SUR strategy can be employed to recover the integer variables with a bounded integer approximation error. We have proved that such an error is related to the MTC linearly. In addition, a combined shrinking and receding horizon strategy has been proposed to satisfy MTC in closed-loop. Recursive feasibility has been proven using a $l$-CI terminal set. An iterative algorithm was developed to explicitly compute a $l$-SRCI set, a specific type of $l$-CI set for switched linear systems. Finally, two numerical examples have been presented to show comparable closed-loop control performance and significant computation speedup using the proposed MPC algorithm, \revv{at the expense of suboptimal solutions}.   

\bibliographystyle{IEEEtran}
\bibliography{main}

\end{document}